\documentclass[A4j,11pt]{article}
\usepackage{graphics}
\usepackage{amsmath}
\usepackage{amssymb}
\usepackage{latexsym}
\usepackage{amsfonts}
\begin{document}

\title{{\bf The constancy of principal curvatures\\
of curvature-adapted submanifolds \\
in symmetric spaces}}
\author{{\bf Naoyuki Koike}}
%
%
\date{}
%
\maketitle
\begin{abstract}
In this paper, we investigate complete curvature-adapted submanifolds with maximal flat 
section and trivial normal holonomy group in symmetric spaces of compact type or 
non-compact type under certain condition, and 
derive the constancy of the principal curvatures of such submanifolds.  
As its result, we can derive that such submanifolds are isoparametric.  
\end{abstract}

\vspace{0.5truecm}

\section{Introduction}
Let $G/K$ be a symmetric space of compact type or non-compact type, and 
$M$ a complete (embedded oriented) Riemannian submanifold in $G/K$.  
Denote by $R$ the curvature tensor of $G/K$.  
Also, denote by $T^{\perp}_xM$ the normal space of $M$ at $x(\in M)$, $A$ the shape tensor 
of $M$, $\nabla^{\perp}$ the normal connection of $M$ and $\exp^{\perp}$ the normal exponential map 
of $M$.  
If, for each $x\in M$ and each $v\in T_x^{\perp}M$, the normal Jacobi operator 
$R(v):=R(\cdot,v)v$ preserves $T_xM$ invariantly and $R(v)\vert_{T_xM}$ commutes with 
$A_v$, then $M$ is said to be {\it curvature-adapted}.  This notion was introduced by 
Berndt-Vanhecke ([3]).  
Curvature-adapted hypersurfaces (in some cases with constant principal curvatures) in rank one symmetric 
spaces were studied by some geometers (see [1,2,5,28] for example).  
If, for each $x\in M$, the normal umbrella $\Sigma_x:=\exp^{\perp}(T^{\perp}_xM)$ 
is totally geodesic, then $M$ is called a {\it submanifold with section}.  
Furthermore, if, for each $x\in M$, the induced metric on $\Sigma_x$ is flat, then $M$ is 
called a {\it submanifold with flat section}.  
Furthermore, if the codimension of $M$ is equal to the rank of $G/K$, then we call 
$M$ a {\it submanifold with maximal flat section}.  
Assume that $M$ is a complete curvature-adapted submanifold with maximal flat section and 
trivial normal holonomy group.  Then, since $M$ has flat section, 
$R(v)\vert_{T_xM}$'s ($v\in T^{\perp}_xM$) commute to one another for each $x\in M$.  
Hence they have the common eigenspace decomposition.  
It is shown that there exist the smooth distributions $D^R_i$ ($i=1,\cdots,m_R$) on $M$ 
such that, for each $x\in M$, $T_xM={\oplus}_{i=1}^{m_R}(D^R_i)_x$ holds and that 
this decomposition is the common eigenspace decomposition of $R(v)\vert_{T_xM}$'s 
($v\in T^{\perp}_xM$).  
Note that $m_R=\mathop{\max}_{v\in T^{\perp}_xM}\sharp\,{\rm Spec}\,R(v)$, 
where $x$ is an arbitrary point of $M$, 
${\rm Spec}(\cdot)$ is the spectrum of $(\cdot)$ and $\sharp(\cdot)$ is the cardinal 
number of $(\cdot)$.  
Let $\pi:\widehat M\to M$ be the universal covering of $M$.  
Then there exist smooth sections $\alpha_i$ of $(\pi^{\ast}T^{\perp}M)^{\ast}$ 
such that, for each $\widehat x\in\widehat M$ and each $v\in T^{\perp}_{\pi(\widehat x)}M$, 
$R(v)\vert_{(D^R_i)_{\pi(\widehat x)}}
=\varepsilon((\alpha_i)_{\widehat x}(v))^2\,{\rm id}$, 
where $\pi^{\ast}T^{\perp}M$ is the induced bundle of 
the normal bundle $T^{\perp}M$ by $\pi$, $(\pi^{\ast}T^{\perp}M)^{\ast}$ is  its 
dual bundle, ${\rm id}$ is the identity transformation of $(D^R_i)_{\pi(\widehat x)}$ and 
$\varepsilon=1$ (resp. $\varepsilon=-1$) in the case where $G/K$ is of compact type 
(resp. of non-compact type).  
Note that each $\alpha_i$ is unique up to the $(\pm 1)$-multiple.  
Set ${\cal R}_M:=\{\pm\alpha_1,\cdots,\pm\alpha_{m_R}\}$ and, for each $x\in M$, define 
a subset ${\cal R}_M^x$ of $(T^{\perp}_xM)^{\ast}$ by 
$${\cal R}_M^x:=\{\pm(\alpha_1)_{\widehat x},\cdots,\pm(\alpha_{m_R})_{\widehat x}\},$$
where $\widehat x$ is an arbitrary point of $\pi^{-1}(x)$.  
Note that ${\cal R}_M^x$ is independent of the choice of $\widehat x\in\pi^{-1}(x)$.  
The system ${\cal R}_M^x$ gives a root system and that it is isomorphic to the 
(restricted) root system of the symmetric pair $(G,K)$.  
Hence, if $\alpha,\beta\in{\cal R}_M$ and if $\beta=F\alpha$ for some $F\in C^{\infty}(M)$, 
then $F=\pm 1\,\,{\rm or}\,\,\pm 2$.  
For convenience, we denote $D^R_i$ by $D^R_{\alpha_i}$.  
For each $\alpha\in{\cal R}_M$ and each parallel normal vector field $\widetilde v$ of $M$, 
we define a function $\alpha(\widetilde v)^2$ over $M$ by 
$$\alpha(\widetilde v)^2(x):=\alpha_{\widehat x}(\widetilde v_x)^2\quad\,\,(x\in M),$$
where $\widehat x$ is an arbitrary point of $\pi^{-1}(x)$.  

On the other hand, since $M$ has flat section and trivial normal holonomy group, 
it follows from the Ricci equation that $A_v$'s ($v\in T^{\perp}_xM$) commute to 
one another for each $x\in M$.  Hence they have the common eigenspace decomposition.  
Set 
$$m_A:=\mathop{\max}_{x\in M}\mathop{\max}_{v\in T^{\perp}_xM}\sharp\,{\rm Spec}\,A_v$$ 
and 
$$U_A:=\{x\in M\,\vert\,\mathop{\max}_{v\in T^{\perp}_xM}\sharp\,{\rm Spec}\,A_v=m_A\}.$$
It is clear that $U_A$ is open in $M$.  Let $U_A^0$ be one of components of $U_A$.  
It is shown that there exist smooth distributions $D^A_i$ ($i=1,\cdots,m_A$) 
on $U_A^0$ such that, for each $x\in U_A^0$, 
$T_xM={\oplus}_{i=1}^{m_A}(D^A_i)_x$ holds and that this decomposition is the common 
eigenspace decomposition of $A_v$'s ($v\in T^{\perp}_xM$).  
Also, there exist smooth sections $\lambda_i$ of $(T^{\perp}M)^{\ast}$ 
($i=1,\cdots,m_A$) such that, 
for each $x\in U_A^0$ and each $v\in T^{\perp}_xM$, 
$A_v\vert_{(D^A_i)_x}=(\lambda_i)_x(v)\,{\rm id}$ holds.  
Set ${\cal A}_M:=\{\lambda_1,\cdots,\lambda_{m_A}\}$.  
For convenience, we denote $D^A_i$ by $D^A_{\lambda_i}$.  

In this paper, we first prove the following results.  

\vspace{0.5truecm}

\noindent
{\bf Theorem A.} {\sl Let $M$ be a complete curvature-adapted submanifold with 
maximal flat section and trivial normal holonomy group in a symmetric space $G/K$ of compact type 
or non-compact type.  Then the following statements {\rm(i)} and {\rm(ii)} hold:

{\rm (i)} For each parallel normal vector field $\widetilde v$ of $M$, 
the eigenvalues $\varepsilon\alpha(\widetilde v)^2$'s ($\alpha\in{\cal R}_M$) of 
$R(\widetilde v)$ are constant over $M$.  

{\rm (ii)} Let $\lambda\in{\cal A}_M$.  Assume that, for any $\alpha\in{\cal R}_M$ with 
$D^R_{\alpha}\cap D^A_{\lambda}\not=\{0\}$, ${\rm dim}(D^R_{\alpha}\cap D^A_{\lambda})\geq 2$.  
Then, for each parallel normal vector field $\widetilde v$ of $U_A^0$, 
the principal curvature $\lambda(\widetilde v)$ of $U_A^0$ for $\widetilde v$ is constant 
along any curve tangent to $D^A_{\lambda}$.}

\vspace{0.5truecm}

In 2006, Heintze-Liu-Olmos ([10]) defined the notion of an {\it isoparametric submanifold} 
in a general Riemannian manifold as a (properly embedded) complete submanifold $M$ with section 
and trivial normal holonomy group satisfying the following condition:

\vspace{0.2truecm}

(Is) Sufficiently close parallel submanifolds of $M$ have constant mean 

curvature with respect to the radial direction.  

\vspace{0.2truecm}

In the sequel, we assume that all isoparametric submanifolds have flat section.   
Next we prove the following result.  

\vspace{0.5truecm}

\noindent
{\bf Theorem B.} {\sl Under the hypothesis of Theorem A, 
assume that, 
for each $\alpha\in{\cal R}_M$, there exists $\lambda\in{\cal A}_M$ such that 
$D^R_{\alpha}\subset D^A_{\lambda}$ holds over $U_A^0$ and that ${\rm dim}\,D^R_{\alpha}\geq 2$.  
Then the following statements {\rm(i)}$-${\rm(iii)} hold:

{\rm(i)} The set $U_A^0$ is equal to $M$ and, for each parallel normal vector field 
$\widetilde v$ of $M$, the principal curvatures $\lambda(\widetilde v)$'s 
($\lambda\in{\cal A}_M$) of $M$ for $\widetilde v$ are constant over $M$.

{\rm(ii)} If $G/K$ is of compact type, then $M$ is isoparametric.  

{\rm(iii)} If $G/K$ is of non-compact type and if $M$ is real analytic, then $M$ is isoparametric.}

\vspace{0.5truecm}

\noindent
{\it Remark 1.1.} 
(i) Principal orbits of a Hermann action $H\curvearrowright G/K$ with 
${\rm cohom}\,H={\rm rank}\,G/K$ are curvature-adapted isoparametric submanifolds 
with maximal flat section and trivial normal holonomy group (see [8,18]).  

(ii) Principal orbits of the isotropy action of any symmetric space of compact type 
(or non-compact type) satisfy the condition that, for each $\alpha\in{\cal R}_M$, there 
exists $\lambda\in{\cal A}_M$ with $D^R_{\alpha}\subset D^A_{\lambda}$ (see [8, Theorem 5.3]).  

(iii) The proofs of Theorem A and Theorem B(i) do not require strong results proved 
in other papers, whereas Theorem B(ii)-(iii) make use of several strong results 
in different papers.  

\vspace{0.5truecm}

By using Theorem B and several strong results in different papers, 
we derive the following result.  

\vspace{0.5truecm}

\noindent
{\bf Theorem C.} {\sl Under the hypothesis of Theorem B, assume that 
$G/K$ is  a simply connected and irreducible symmetric space of compact type and rank 
greater than one.  Then $M$ is congruent to a principal orbit of the isotropy action 
of $G/K$.}

\vspace{0.5truecm}

\noindent
{\bf Theorem D.} {\sl Under the hypothesis of Theorem B, assume that 
$G/K$ is an irreducible symmetric space of non-compact type and rank greater than one, 
and that $M$ is real anlaytic and has no non-Euclidean type focal point on the ideal 
boundary of $G/K$.  Then $M$ is a principal orbit of a Hermann action on $G/K$.}

\vspace{0.5truecm}

The notion of a non-Euclidean type focal point on the ideal boundary of $G/K$ was introduced in [19].  
See the next section about the definition of this notion.  

\vspace{0.5truecm}

\noindent
{\it Remark 1.2.} 
(i) The principal orbits of any Hermann action on a symmetric space $G/K$ of 
non-compact type have no non-Euclidean type focal point on the ideal boundary of $G/K$.  

(ii) Let $G/K$ be an irreducible symmetric space of non-compact type and rank greater 
than one, $N$ the nilpotent part in the Iwasawa's decomposition $G=KAN$ of $G$ and 
$M$ a principal orbit of the $N$-action on $G/K$.  
Assume that the multiplicity of each root of the (restricted) root system of the 
symmetric pair $(G,K)$ is greater than one.  Then $M$ satisfies all the hypothesis of 
Theorem B (see [20]) and it is real analytic.  However it has a non-Euclidean type focal 
point on the ideal boundary of $G/K$.  
On the other hand, we can show that $M$ does not occur as a principal orbit of a Hermann 
action on $G/K$.  Thus, in Theorem D, is indispensable the condition that $M$ has no 
non-Euclidean type focal point on the ideal boundary of $G/K$.  

(iii) There exists a Hermann action on an irreducible symmetric space $G/K$ of non-compact type and 
rank greater than one such that its principal orbits satisfy all the hypothesis in Theorem D 
but that it is not conjugate to the isotropy actions of $G/K$.  
For example, see Table 1 in Section 4 about all of such Hermann actions 
on irreducible rank two symmetric spaces of non-compact type.  

\section{Basic notions and facts}
In this section, we shall recall the basic notions and facts in a symmetric space.  
See [7, Pages 94-95] or [30, Page 177] about the alebraic structure of a symmetric space 
and the Jacobi field on the space.  
We use the notations in Introduction.  
Let $M$ be an embedded submanifold in a Riemannian manifold $N$.
Denote by $\nabla,\,\widetilde{\nabla}$ and $A$ the Riemannian connection of $M$, that 
of $N$ and the shape tensor of $M$, respectively.  
Take a unit normal vector $v$ of $M$ at $x$ and denote by $\gamma_v$
the geodesic in $N$ with $\gamma'_v(0)=v$, where $\gamma'_v(0)$ is the velocity vector of 
$\gamma_v$ at $0$.  
If there exists an $M$-Jacobi field $J$ along $\gamma_v$ satisfying
$J(0)\not=0$ and $J(s_0)=0$, then the real number $s_0$ is called a
{\it focal radius along} $\gamma_v$.
We consider the case where $N$ is a symmetric space $G/K$.
In this case, 
the strongly $M$-Jacobi field $J$ along $\gamma_v$ with 
$J(0)=X$ (hence 
$\frac{\widetilde{\nabla}}{\partial s}J\vert_{s=0}=-A_vX$) 
is given by 
$$J(s)=\left(P_{\gamma_v\vert_{[0,s]}}\circ\left(\cos(s\sqrt{R(v)})
-\frac{\sin(s\sqrt{R(v)})}{\sqrt{R(v)}}\circ A_v\right)\right)(X),\leqno{(2.1)}$$
where $\frac{\widetilde{\nabla}}{\partial s}$ is the covariant derivative 
along $\gamma_v$ with respect to $\widetilde{\nabla}$ and $P_{\gamma_v\vert_{[0,s]}}$ is 
the parallel translation along $\gamma_v\vert_{[0,s]}$.  
In the case where $M$ has flat section, any focal radius of $M$ along $\gamma_v$ is given 
as a zero point of a strongly $M$-Jacobi field along $\gamma_v$.  
Hence the set of all focal radii of $M$ along $\gamma_v$ coincides with the zero point set of 
the real-valued function $F_v$ over ${\Bbb R}$ defined by 
$$F_v(s):={\rm det}\left(\cos(s\sqrt{R(v)})
-\frac{\sin(s\sqrt{R(v)})}{\sqrt{R(v)}}\circ A_v\right).$$
In 1995, Terng-Thorbergsson ([32]) defined the notion of an {\it equifocal submanifold} as a 
compact submanifold with flat section and trivial normal holonomy group satisfying 
the following condition:

\vspace{0.2truecm}

(PF) $M$ has parallel focal structure, that is,
for any parallel normal vector 

field $\widetilde v$ of $M$, the focal radii along $\gamma_{\widetilde v_x}$ 
are independent of the choice of 

$x\in M$ (with considering their multiplicities).

\vspace{0.2truecm}

\noindent
Let $H$ be a closed subgroup of $G$.  The $H$-action on $G/K$ is called a
{\it polar action} if $H$ is compact and if there exists a complete embedded
submanifold $\Sigma$ meeting all principal $H$-orbits orthogonally.
The submanifold $\Sigma$ is called a {\it section} of this action.  
Furthermore, if the induced metric on $\Sigma$ is flat, then the $H$-action
is called a {\it hyperpolar action}.  
It is known that principal orbits of a hyperpolar action are equifocal.  
On the other hand, in 1995, E. Heintze, R.S. Palais, C.L. Terng and
G. Thorbergsson ([11]) proved that, in the case where $G/K$ is a simply connected symmetric space 
of compact type, any homogeneous submanifold with flat section in $G/K$ is a principal orbit of 
a hyperpolar action.  If $G/K$ is of compact type and if there exists an involution $\sigma$ of 
$G$ with $({\rm Fix}\,\sigma)_0\subset H\subset{\rm Fix}\,\sigma$, then the 
$H$-action on $G/K$ is called a {\it Hermann action}.  It is known that
Hermann actions are hyperpolar ([11]) and that the principal orbits of a Hermann action are 
curvature-adapted (see [8]).  
In 2001, A. Kollross ([24]) proved that, in the case where $G/K$ is an irreducible simply connected 
symmetric space of compact type, 
hyperpolar actions of cohomogeneity greater than one on an irreducible simply
connected symmetric space of compact type are orbit equivalent to Hermann actions.  
In 2002, U. Christ ([6]) proved that, in the case where $G/K$ is an irreducible simply connected 
symmetric space of compact type, any irreducible equifocal submanifold of codimension greater than 
one in $G/K$ is homogeneous.  
Note that there was a gap in his proof but, in 2012, C. Gorodski and E. Heintze ([9]) 
closed the gap.  

Therefore we obtain the following fact.

\vspace{0.5truecm}

\noindent
{\bf Fact 2.1.} {\sl Any equifocal submanifold of codimension greater than one 
in any irreducible simply connected symmetric space of compact type is a principal orbit of 
a Hermann action.}

\vspace{0.5truecm}

Heintze-Liu-Olmos ([10]) showed the following fact.  

\vspace{0.5truecm}

\noindent
{\bf Fact 2.2.} {\sl For a compact submanifold with flat section in a symmetric space of compact 
type, it is equifocal if and only if it is isoparametric.}

\vspace{0.5truecm}

\noindent
{\it Remark 2.1.} In more general, it is shown that, 
for a (not necessarily compact) submanifold with flat section and trivial normal holonomy group 
in a symmetric space of compact type, it has parallel focal structure if and only if it is isoparametric.

\vspace{0.5truecm}

When a non-compact submanifold $M$ in a symmetric space $G/K$ of non-compact type deforms 
as its principal curvatures approach to zero, 
its focal set vanishes beyond the ideal boundary $(G/K)(\infty)$ of $G/K$.  
For example, when an open potion of a totally umbilic sphere in a hyperbolic space of constant 
curvature $c(<0)$ deforms as its principal curvatures approach to $\sqrt{-c}$, its focal point 
approach to $(G/K)(\infty)$ and, when it furthermore deforms as its principal curvatures approach to 
a positive value smaller than $\sqrt{-c}$, the focal point vanishes beyond $(G/K)(\infty)$.  
On the base of this fact, we recognized that, for a non-compact submanifold 
in a symmetric space of non-compact type, the parallelity of the focal 
structure is not an essential condition.  
So, in 2004, we [16] introduced the notion of a complex focal radius along 
the normal geodesic $\gamma_v$ of such a submanifold as a general notion of 
a focal radius along $\gamma_v$.  This notion was defined as 
the zero points of the complex-valued function $F_v^{\Bbb C}$ over ${\Bbb C}$ 
defined by 
$$F_v^{\Bbb C}(z):={\rm det}\left(\cos(z\sqrt{R(v)^{\Bbb C}})
-\frac{\sin(z\sqrt{R(v)^{\Bbb C}})}{\sqrt{R(v)^{\Bbb C}}}\circ A_v^{\Bbb C}
\right),$$
where $R(v)^{\Bbb C}$ (resp. $A_v^{\Bbb C}$) is the complexification of $R(v)$ 
(resp. $A_v$).  
In the case where $M$ is real analytic, we [17] showed that
complex focal radii along $\gamma_v$ indicate the positions of focal points of
the extrinsic complexification $M^{\Bbb C}(\hookrightarrow G^{\Bbb C}/K^{\Bbb C})$ of $M$ 
along the complexified geodesic $\gamma_{\iota_{\ast}v}^{\Bbb C}$, where
$G^{\Bbb C}/K^{\Bbb C}$ is the anti-Kaehlerian symmetric space associated with
$G/K$ and $\iota$ is the natural immersion of $G/K$ into
$G^{\Bbb C}/K^{\Bbb C}$.  

We ([16]) defined the notion of a {\it complex equifocal submanifold} as a 
(properly embedded) complete submanifold with flat section and trivial normal holonomy 
group satisfying the following condition:

\vspace{0.2truecm}

(PCF) $M$ has parallel complex focal structure, that is, for any parallel 

normal vector field $\widetilde v$ of $M$, the complex focal radii along
$\gamma_{\widetilde v_x}$ are indepe-

ndent of the choice of $x\in M$ (considering their multiplicities).

\vspace{0.2truecm}

\noindent
We should call this submanifold a {\it equi-complex-focal submanifold} but 
called a {\it complex equifocal submanifold} for simplicity.  
We ([17]) showed the following fact.  

\vspace{0.5truecm}

\noindent
{\bf Fact 2.3.} {\sl Let $M$ be a complete submanifold with flat section in a symmetric space 
$G/K$ of non-compact type.  If $M$ is isoparametric, then it is complex equifocal.  Convesely, if 
$M$ is real analytic, complex equifocal and curvature-adapted, then it is isoparametric.}

\vspace{0.5truecm}

Let $G/K$ be a symmetric space of non-compact type and $H$ a closed subgroup of $G$.  
The $H$-action on $G/K$ is called a {\it polar action} 
if there exists a complete embedded submanifold $\Sigma$ meeting all principal $H$-orbits 
orthogonally.  Furthermore, if the induced metric on $\Sigma$ is flat, then the $H$-action 
is called a {\it hyperpolar action}.  Note that a polar action (resp. a hyperpolar action) 
on a symmetric space of non-compact type was called a {\it complex polar action} (resp. 
{\it complex hyperpolar action}) in [17,18,20].  
In [17], it was proved that principal orbits of a hyperpolar action on $G/K$ are complex 
equifocal and that any homogeneous submanifold with flat section in $G/K$ is a principal orbit of 
a hyperpolar action.  If there exists an involution $\sigma$ of
$G$ with $({\rm Fix}\,\sigma)_0\subset H\subset{\rm Fix}\,\sigma$, then the
$H$-action on $G/K$ is called a {\it Hermann type action}.  
For simplicity, we call this action a {\it Hermann action} in this paper.  
It is easy to show that Hermann actions are hyperpolar.  

At the end of this section, we recall the notion of a non-Euclidean type focal point on 
the ideal boundary of a Hadamard manifold introduced in [19].  
Let $N$ be a Hadamard manifold, $N(\infty)$ the ideal boundary of $N$ and $M$ 
a submanifold in $N$.  Take $v\in T^{\perp}_xM$.  
Let $\gamma_v:[0,\infty)\to N$ be the normal geodesic of $M$ of direction $v$.  
If there exists an $M$-Jacobi field $J$ along $\gamma_v$ satisfying 
$\lim\limits_{s\to\infty}\frac{\vert\vert J(s)\vert\vert}{s}=0$, then we call 
$\gamma_v(\infty)\,(\in N(\infty))$ a {\it focal point} of $M$ {\it on the ideal boundary} 
$N(\infty)$ {\it along} $\gamma_v$, where $\gamma_v(\infty)$ is the asymptotic class of 
$\gamma_v$.  In fact, if such an $M$-Jacobi field $J$ exists, then for the standard geodesic 
variation $\delta:[0,\infty)\times(-\varepsilon,\varepsilon)\to N$ having $J$ as the variational 
vector field, the asymptotic classes $\gamma_s(\infty)$'s of 
$\gamma_s:t\mapsto\delta(t,s)$ $(-\varepsilon<s<\varepsilon)$ coincide with $\gamma_v(\infty)$, 
where $\varepsilon$ is a positive number.  
Hence $\gamma_v(\infty)$ should be interpreted as a focal point of $M$.  
Also, if there exists an $M$-Jacobi field $J$ along $\gamma_v$ satisfying 
$\lim\limits_{s\to\infty}\frac{\vert\vert J(s)\vert\vert}{s}=0$ and 
${\rm Sec}(v,J(0))\not=0$, then we call $\gamma_v(\infty)$ a 
{\it non-Euclidean type focal point} $M$ {\it on} $N(\infty)$ 
{\it along} $\gamma_v$, where ${\rm Sec}(v,J(0))$ is the 
sectional curvature for the $2$-plane spanned by $v$ and $J(0)$.  
If, for any nonzero normal vector $v$ of $M$, 
$\gamma_v(\infty)$ is not a non-Euclidean type focal point of $M$ on $N(\infty)$ 
along $\gamma_v$, then we say that $M$ 
{\it has no non-Euclidean type focal point on the ideal boundary}.  

\section{Proofs of Theorems A, B and C} 
Let $M$ be a complete curvature-adapted submanifold in $G/K$ as in the statement of 
Theorem A.  We shall use the notations in Introduction.  Denote by $\nabla$ and 
$\widetilde{\nabla}$ the Riemannian connections of $M$ and $G/K$, respectively.  
In the sequel, for each $\alpha\in{\cal R}_M$ with 
$2\alpha\notin{\cal R}_M$, $D^R_{2\alpha}$ implies 
the zero distribution.  
First we note that the following relations hold:
$$\begin{array}{c}
\displaystyle{R(D^R_{\alpha},T^{\perp}M)T^{\perp}M\subset D^R_{\alpha},\quad\,\,
R(D^R_{\alpha},D^R_{\beta})T^{\perp}M\subset D^R_{\alpha+\beta}\oplus D^R_{\alpha-\beta},}\\
\displaystyle{R(D^R_{\alpha},D^R_{\alpha})T^{\perp}M\subset D^R_{2\alpha},}\\
\displaystyle{R(D^R_{\alpha},T^{\perp}M)D^R_{\beta}\subset D^R_{\alpha+\beta}\oplus D^R_{\alpha-\beta}
\oplus T^{\perp}M.}
\end{array}\leqno{(3.1)}$$
because $M$ has maximal flat section.  

We shall prove Theorem A.  

\vspace{0.5truecm}

\noindent
{\it Proof of Theorem A.} 
First we shall show the statement (i).  Take a parallel normal vector field $\widetilde v$ of $M$.  
Let $V$ be a sufficiently small open set of $M$.  
Take a local unit section $\widetilde X$ of $TM\ominus D^R_{2\alpha}$ defined over $V$ and 
a local unit section $\widetilde Y$ of $D^R_{\alpha}$ defined over $V$, where 
$TM\ominus D^R_{2\alpha}$ means $TM\cap(D^R_{2\alpha})^{\perp}$.  
Then we have 
$$\widetilde{\nabla}_{\widetilde X}
\left(R(\widetilde Y,\widetilde v)\widetilde v\right)
=\varepsilon\widetilde X(\alpha(\widetilde v)^2)\widetilde Y
+\varepsilon\alpha(\widetilde v)^2
\widetilde{\nabla}_{\widetilde X}\widetilde Y.\leqno{(3.2)}$$
On the other hand, we have 
$$
\widetilde{\nabla}_{\widetilde X}
\left(R(\widetilde Y,\widetilde v)\widetilde v\right)
=R(\nabla_{\widetilde X}\widetilde Y,\widetilde v)\widetilde v
-R(\widetilde Y,A_{\widetilde v}\widetilde X)\widetilde v
-R(\widetilde Y,\widetilde v)A_{\widetilde v}\widetilde X,\leqno{(3.3)}$$
where we use $R(\widetilde v)\vert_{T^{\perp}M}=0$.  
From $\langle\nabla_{\widetilde X}\widetilde Y,\widetilde Y\rangle=0$, we have 
$\langle R(\nabla_{\widetilde X}\widetilde Y,\widetilde v)\widetilde v,\widetilde Y
\rangle=0$.  
Hence, by taking the inner product of $(3.2)$ and $(3.3)$ with $\widetilde Y$, we obtain 
$$\varepsilon\widetilde X(\alpha(\widetilde v)^2)
=-\langle R(\widetilde Y,A_{\widetilde v}\widetilde X)\widetilde v,\widetilde Y\rangle
-\langle R(\widetilde Y,\widetilde v)A_{\widetilde v}\widetilde X,\widetilde Y\rangle.
\leqno{(3.4)}$$
Also, since $M$ is curvature-adapted, $A_{\widetilde v}\widetilde X$ is a local section of 
$TM\ominus D^R_{2\alpha}$.  Hence 
it follows from $(3.1)$ that 
$$R(\widetilde Y,A_{\widetilde v}\widetilde X)\widetilde v\in(D^R_{\alpha})^{\perp}\,\,\,\,
{\rm and}\,\,\,\,R(\widetilde Y,\widetilde v)A_{\widetilde v}\widetilde X
\in(D^R_{\alpha})^{\perp}.
\leqno{(3.5)}$$
From $(3.4)$ and $(3.5)$, we obtain 
$$\widetilde X(\alpha(\widetilde v)^2)=0.\leqno{(3.6)}$$
Take local unit sections $\widetilde Z_i$ ($i=1,2$) of $D^R_{2\alpha}$ defined over $V$.  
Then, in similar to $(3.4)$, we can show 
$$\varepsilon\widetilde Z_1(2\alpha(\widetilde v)^2)
=-\langle R(\widetilde Z_2,A_{\widetilde v}\widetilde Z_1)\widetilde v,\widetilde Z_2\rangle
-\langle R(\widetilde Z_2,\widetilde v)A_{\widetilde v}\widetilde Z_1,\widetilde Z_2\rangle.
\leqno{(3.7)}$$
Also, it follows from $(3.1)$ that 
$$R(\widetilde Z_2,A_{\widetilde v}\widetilde Z_1)\widetilde v=0\quad\,\,{\rm and}\quad\,\,
R(\widetilde Z_2,\widetilde v)A_{\widetilde v}\widetilde Z_1\in T^{\perp}M.$$
Hence we obtain 
$$\widetilde Z_1(\alpha(\widetilde v)^2)=0.\leqno{(3.8)}$$
From $(3.6),\,(3.8)$, the arbitrarinesses of $\widetilde X, \widetilde Z_1$ and $V$, 
it follows that $\alpha(\widetilde v)^2$ is constant over $M$.  

Next we shall show the statement (ii).  
Take $x\in U_A^0$ and $\alpha\in{\cal R}_M$ with 
$D^R_{\alpha}\cap D^A_{\lambda}\not=\{0\}$.  
Take local sections $\widetilde X_i$ ($i=1,2$) of 
$D^R_{\alpha}\cap D^A_{\lambda}$ over a neighborhood $V$ of $x$ in $U_A^0$ 
which give an orthonormal system at each point of $V$, where 
we use ${\rm dim}(D^R_{\alpha}\cap D^A_{\lambda})\geq 2$.  
Then we have 
$$\begin{array}{l}
\displaystyle{R((\widetilde X_1)_x,(\widetilde X_2)_x)\widetilde v_x
=\widetilde{\nabla}_{(\widetilde X_1)_x}\widetilde{\nabla}_{\widetilde X_2}\widetilde v
-\widetilde{\nabla}_{(\widetilde X_2)_x}\widetilde{\nabla}_{\widetilde X_1}\widetilde v
-\widetilde{\nabla}_{[\widetilde X_1,\widetilde X_2]_x}\widetilde v}\\
\hspace{3.35truecm}\displaystyle{=-((\widetilde X_1)_x(\lambda(\widetilde v)))
(\widetilde X_2)_x
+((\widetilde X_2)_x(\lambda(\widetilde v)))(\widetilde X_1)_x}\\
\hspace{3.85truecm}\displaystyle{
-(\lambda)_x(\widetilde v_x)[\widetilde X_1,\widetilde X_2]_x
+A_{\widetilde v_x}([\widetilde X_1,\widetilde X_2]_x).}
\end{array}\leqno{(3.9)}$$
On the other hand, since $(\widetilde X_i)_x$ ($i=1,2$) belong to $(D^R_{\alpha})_x$ 
and $\widetilde v_x$ belongs to $T^{\perp}_xM$, it follows from $(3.1)$ that 
$$R((\widetilde X_1)_x,(\widetilde X_2)_x)\widetilde v_x\in(D^R_{2\alpha})_x.
\leqno{(3.10)}$$
From $(3.9)$ and $(3.10)$, we have 
$$\begin{array}{l}
\displaystyle{
((\widetilde X_1)_x(\lambda(\widetilde v)))(\widetilde X_2)_x
-((\widetilde X_2)_x(\lambda(\widetilde v)))(\widetilde X_1)_x}\\
\displaystyle{\equiv\left(A_{\widetilde v_x}-(\lambda)_x(\widetilde v_x)\right)
([\widetilde X_1,\widetilde X_2]_x)\quad\,\,({\rm mod}\,\,(D^R_{2\alpha})_x).}
\end{array}\leqno{(3.11)}$$
The left-hand side of this relation belongs to $(D^A_{\lambda})_x$ but 
the right-hand side of this relation is orthogonal to $(D^A_{\lambda})_x$.  
Hence, it follows that the left-hand side of $(3.11)$ vanishes.  
Therefore, it follows from the linear independency of $(\widetilde X_1)_x$ and 
$(\widetilde X_2)_x$ that $(\widetilde X_1)_x(\lambda(\widetilde v))=0$.  
From the arbitrarinesses of $x$ and $\widetilde X_1$, it follows that 
$\lambda(\widetilde v)$ is constant along any curve tangent to 
$D^R_{\alpha}\cap D^A_{\lambda}$.  Furthermore, the statement (ii) of Theorem A follows 
from the arbitrariness of $\alpha$, where we use the curvature-adaptedness of $M$.  
\hspace{1.5truecm}q.e.d.

\vspace{0.3truecm}

Next we shall prove the statement (i) of Theorem B.  

\vspace{0.3truecm}

\noindent
{\it Proof of (i) of Theorem B.} 
Let $D^R_{\alpha}\subset D^A_{\lambda}$ and set 
$D^R_{\alpha,2\alpha}:=D^R_{\alpha}\oplus D^R_{2\alpha}$.  

(Step I) First we shall show that $D^R_{\alpha,2\alpha}$ is a totally geodesic distribution 
on $U_A^0$.  Take a parallel normal vector field $\widetilde v$ of $U_A^0$.  
Fix $x_0\in U_A^0$.  Take local sections $\widetilde X_1$ and $\widetilde X_2$ of 
$D^R_{\alpha}$ defined over a neighborhood $V$ of $x_0$ in $U_A^0$.  
Since $\alpha(\widetilde v)^2$ is constant over $M$ by the statement (i) of Theorem A, 
we have 
$$\widetilde{\nabla}_{(\widetilde X_1)_x}(R(\widetilde X_2,\widetilde v)\widetilde v)
=\varepsilon\alpha(\widetilde v)^2(x)
\widetilde{\nabla}_{(\widetilde X_1)_x}\widetilde X_2
\leqno{(3.12)}$$
On the other hand, since $M$ has flat section, we obtain 
$$\begin{array}{l}
\displaystyle{\widetilde{\nabla}_{(\widetilde X_1)_x}(R(\widetilde X_2,\widetilde v)
\widetilde v)
=R(\nabla_{(\widetilde X_1)_x}\widetilde X_2,\widetilde v_x)\widetilde v_x}\\
\hspace{3.82truecm}\displaystyle{
-\lambda_x(\widetilde v_x)
R((\widetilde X_2)_x,(\widetilde X_1)_x)\widetilde v_x}\\
\hspace{3.82truecm}
\displaystyle{-\lambda_x(\widetilde v_x)
R((\widetilde X_2)_x,\widetilde v_x)(\widetilde X_1)_x.}
\end{array}\leqno{(3.13)}$$
It follows from $(3.1)$ that 
$$R((\widetilde X_2)_x,(\widetilde X_1)_x)\widetilde v_x\in(D^R_{2\alpha})_x\leqno{(3.14)}$$
and
$$R((\widetilde X_2)_x,\widetilde v_x)(\widetilde X_1)_x\in 
T^{\perp}_xM\oplus(D^R_{2\alpha})_x.
\leqno{(3.15)}$$
Also, since $M$ has flat section, we have 
$$R(\nabla_{(\widetilde X_1)_x}\widetilde X_2,\widetilde v_x)\widetilde v_x\in T_xM.
\leqno{(3.16)}$$
From $(3.12)-(3.16)$, we obtain 
$$R(\nabla_{(\widetilde X_1)_x}\widetilde X_2,\widetilde v_x)\widetilde v_x
\equiv\varepsilon\alpha(\widetilde v)^2(x)
\nabla_{(\widetilde X_1)_x}\widetilde X_2\quad\,\,({\rm mod}\,(D^R_{2\alpha})_x).
\leqno{(3.17)}$$
Therefore, from the arbitrariness of $\widetilde v$, we obtain 
$\nabla_{(\widetilde X_1)_x}\widetilde X_2\in(D^R_{\alpha,2\alpha})_x$.  
Furthermore, from the arbitarinesses of $\widetilde X_1,\widetilde X_2$ and $x$, 
it follows that $\nabla_{D^R_{\alpha}}D^R_{\alpha}\subset 
D^R_{\alpha,2\alpha}$ holds on $V$.  Furthermore, it follows from the arbitrariness of 
$x_0$ that $\nabla_{D^R_{\alpha}}D^R_{\alpha}\subset D^R_{\alpha,2\alpha}$ holds 
on $U_A^0$.  Similarly we can show that 
$\nabla_{D^R_{2\alpha}}D^R_{2\alpha}\subset D^R_{2\alpha}$, 
$\nabla_{D^R_{\alpha}}D^R_{2\alpha}\subset D^R_{\alpha,2\alpha}$ and 
$\nabla_{D^R_{2\alpha}}D^R_{\alpha}\subset D^R_{\alpha,2\alpha}$ hold on $U_A^0$.  
Therefore $D^R_{\alpha,2\alpha}$ is a totally geodesic distribution on $U_A^0$.  

(Step II) 
Since $D^R_{\alpha,2\alpha}\vert_{U_A^0}$ is totally geodesic, it is integrable.  
Denote by $\mathfrak F_{\alpha,2\alpha}$ the foliation on $U_A^0$ 
whose leaves are the integral manifolds of $D^R_{\alpha,2\alpha}\vert_{U_A^0}$ and 
$L^{\alpha,2\alpha}_x$ the leaf of $\mathfrak F_{\alpha,2\alpha}$ through $x(\in U_A^0)$.  
Let $\xi_{\alpha_{\widehat x}}$ be the element of $T^{\perp}_xM$ defined by 
$\langle \xi_{\alpha_{\widehat x}},\,\cdot\rangle=\alpha_{\widehat x}(\cdot)$.  
Set $E^x_{\alpha}:=(D^R_{\alpha,2\alpha})_x\oplus{\rm Span}\{\xi_{\alpha_{\widehat x}}\}$ and 
$\Sigma^{\alpha}_x:=\exp_x(E^x_{\alpha})$, where $\exp_x$ is the exponential map of $G/K$ at $x$.  
In this step, we shall show that $L^{\alpha,2\alpha}_x=\Sigma^{\alpha}_x\cap U_A^0$.  
Take $\widehat x\in\pi^{-1}(x)$.  
Let $K_x$ be the isotropy subgroup of $G$ at $x$ and $\theta_x$ an involution of $G$ with 
$({\rm Fix}\,\theta_x)_0\subset K_x\subset{\rm Fix}\,\theta_x$, where 
${\rm Fix}\,\theta_x$ is the fixed point group of $\theta_x$ and 
$({\rm Fix}\,\theta_x)_0$ is its identity component.  
Denote by the same symbol $\theta_x$ the involution of $\mathfrak g$ induced from 
$\theta_x$.  Let $\mathfrak g=\mathfrak k_x+\mathfrak p_x$ be the eigenspace decomposition 
of $\theta_x$.  
The subspace $\mathfrak p_x$ is identified with $T_x(G/K)$.  It is easy to show that 
$R(E^x_{\alpha},E^x_{\alpha})E^x_{\alpha}\subset E^x_{\alpha}$, which 
implies that $E_x$ is a Lie triple system of $\mathfrak g$.  
Therefore $\Sigma^{\alpha}_x$ is totally geodesic in $G/K$ (see P. 224 of [12]).  
In case of $2\alpha\in{\cal R}_M$, 
$\Sigma^{\alpha}_x$ is a simply connected rank one symmetric space of two sectional 
curvatures $\varepsilon\alpha_{\widehat x}(\xi_{\alpha_{\widehat x}})$ and 
$4\varepsilon\alpha_{\widehat x}(\xi_{\alpha_{\widehat x}})$.  
In case of $2\alpha\notin{\cal R}_M$, 
$\Sigma^{\alpha}_x$ is a real space form of constant curvature 
$\varepsilon\alpha_{\widehat x}(\xi_{\alpha_{\widehat x}})$.  
It is clear that the tangent spaces of $L^{\alpha,2\alpha}_x$ and 
$\Sigma^{\alpha}_x\cap U_A^0$ at $x$ coincide with each other.  Since 
$D^R_{\alpha,2\alpha}$ is totally geodesic, $L^{\alpha,2\alpha}_x$ is a totally geodesic 
submanifold in $U_A^0$.  Let $\widetilde X_1$ and $\widetilde X_2$ be local tangent vector 
fields of $\Sigma^{\alpha}_x\cap U_A^0$ defined over a neighborhood $V$ of $x$ in 
$\Sigma^{\alpha}_x\cap U_A^0$.  
Since $\Sigma^{\alpha}_x$ is totally geodesic in $G/K$, we have 
$\widetilde{\nabla}_{(\widetilde X_1)_x}\widetilde X_2$ is tangent to $\Sigma^{\alpha}_x$.  
Hence $\nabla_{(\widetilde X_1)_x}\widetilde X_2$ is tangent to 
$\Sigma^{\alpha}_x\cap U_A^0$.  This implies that $\Sigma^{\alpha}_x\cap U_A^0$ is totally 
geodesic in $M$.  From these facts, we can derive 
$L^{\alpha,2\alpha}_x=\Sigma^{\alpha}_x\cap U_A^0$.  

(Step III) Take any two points $x$ and $y$ of $U_A^0$.  
Since $L^{\alpha,2\alpha}_x=\Sigma^{\alpha}_x\cap U_A^0$ and 
$L^{\alpha,2\alpha}_y=\Sigma^{\alpha}_y\cap U_A^0$, 
$L^{\alpha,2\alpha}_x$ (resp. $L^{\alpha,2\alpha}_y$) is 
a hypersurface in $\Sigma^{\alpha}_x$ (resp. $\Sigma^{\alpha}_y$).  
We shall compare these hypersurfaces.  
Let $\xi_{\alpha}$ be a unit normal vector field of $\widehat M$ such that 
$(\xi_{\alpha})_{\widehat z}=\xi_{\alpha_{\widehat z}}$ for any $\widehat z\in\widehat M$.  
According to (i) of Theorem A, $\alpha$ is parallel with respect to the normal connection 
of $\widehat M$.  Hence so is also $\xi_{\alpha}$.  
According to the assumption, $D^R_{2\alpha}\subset D^A_{\bar{\lambda}}$ hold 
over $U_A^0$ for some $\bar{\lambda}\in{\cal A}_M$.  
Since $\Sigma^{\alpha}_x$ and $\Sigma^{\alpha}_y$ are totally geodesic in $G/K$, 
it follows from $D^R_{\alpha}\subset D^A_{\lambda}$ and $D^R_{2\alpha}\subset 
D^A_{\bar{\lambda}}$ that $\lambda(\widetilde w)=\bar{\lambda}(\widetilde w)=0$ 
for any parallel normal vector field $\widetilde w$ of $M$ orthogonal to $\xi_{\alpha}$.  
On the other hand, since $\Sigma^{\alpha}_x$ is totally geodesic in $G/K$, 
the shape operator of the hypersurface $L^{\alpha,2\alpha}_x$ in $\Sigma^{\alpha}_x$ 
coincides with $A_{\xi_{\alpha}}\vert_{TL^{\alpha,2\alpha}_x}$.  Hence 
$L^{\alpha,2\alpha}_x$ is a curvature-adapted hypersurface in $\Sigma_x^{\alpha}$.  
Also, it is clear that the hypersurface $L^{\alpha,2\alpha}_x$ has at most two 
principal curvatures $\lambda(\xi_{\alpha})$ and $\bar{\lambda}(\xi_{\alpha})$, and, 
furthermore, it follows from (ii) of Theorem A that they are constant along leaves of 
the corresponding eigendistributions.  
Since $\Sigma^{\alpha}_x$ is a rank one symmetric space and $L^{\alpha,2\alpha}_x$ is a complete 
curvature adapted hypersurface with two distinct principal curvatures satisfying a special condition that 
the eigendistributions of the normal Jacobi operator coincide with those of the shape operator, 
this hypersurface $L^{\alpha,2\alpha}_x$ has constant principal curvatures (see [14,15,26,27]).  
Hence $\lambda(\xi_{\alpha})$ and $\bar{\lambda}(\xi_{\alpha})$ are constant along 
$L^{\alpha,2\alpha}_x$.  
Similarly, $\lambda(\xi_{\alpha})$ and $\bar{\lambda}(\xi_{\alpha})$ are constant along 
$L^{\alpha,2\alpha}_y$.  
Since $\alpha(\xi_{\alpha})^2$ is constant over $M$, it follows that $\Sigma^{\alpha}_x$ 
and $\Sigma^{\alpha}_y$ are isometric to each other.  
Also, since ${\cal F}_{\alpha,2\alpha}$ is a totally geodesic foliation, it follows from the result 
in [4] that the element of holonomy along any curve orthogonal to leaves of ${\cal F}_{\alpha,2\alpha}$ 
consists of local isometries between the leaves.  See [4] about the definition of 
an element of holonomy.  Furthermore, since $M$ is complete, 
the orthogonal complementary distribution of ${\cal F}_{\alpha,2\alpha}$ is an Ehresmann connection 
for ${\cal F}_{\alpha,2\alpha}$ in the sense of [4] and hence any two leaves of 
${\cal F}_{\alpha,2\alpha}$ is joined by a curve orthogonal to the leaves.  Therefore 
it follows that $L^{\alpha,2\alpha}_x$ and $L^{\alpha,2\alpha}_y$ are locally isometric to each other.  
From these facts, we can derive $\lambda_x((\xi_{\alpha})_x)=\lambda_y((\xi_{\alpha})_y)$.  
From the arbitrarinesses of $x$ and $y$, $\lambda(\xi_{\alpha})$ is constant over $U_A^0$.  
After all $\lambda(\widetilde v)$ is constant over $U_A^0$ for any parallel normal vector 
field $\widetilde v$ of $M$.  
Furthermore, from this fact, $U_A^0=M$ is derived directly.  
This completes the proof.  
\hspace{1.5truecm}q.e.d.

\vspace{0.5truecm}

Next we shall prove the statements (ii) and (iii) of Theorem B.  

\vspace{0.5truecm}

\noindent
{\it Proof of (ii)  and (iii) of Theorem B.} 
The focal radii and the complex focal radii along the normal geodesic $\gamma_v$ of $M$ 
are given as the zero points of the functions $F_v$ and $F_v^{\Bbb C}$ as in Section 2.  
Hence, it follows from the additional assumption in Theorem B that, in the case where 
$G/K$ is of compact type, the set of all the focal radii along $\gamma_v$ is equal to 
$$
\left\{\left.\frac{{\rm arctan}(\alpha_{\hat x}(v)/\lambda_x(v))+j\pi}{\alpha_{\hat x}(v)}
\,\right\vert\,(\alpha,\lambda)\in{\cal R}_M\times{\cal A}_M\,\,{\rm s.t.}\,\,
D_{\alpha}\subset D_{\lambda},\,\,\,j\in{\Bbb Z}\right\}$$
and, in the case where $G/K$ is of non-compact type, the set of all the complex focal radii along 
$\gamma_v$ is equal to 
$$\begin{array}{l}
\displaystyle{
\left\{\left.
\frac{{\rm arctanh}(\alpha_{\hat x}(v)/\lambda_x(v))+j\pi\sqrt{-1}}{\alpha_{\hat x}(v)}
\,\right\vert\,(\alpha,\lambda)\in{\cal R}_M\times{\cal A}_M\,\,{\rm s.t.}\,\,\right.}\\
\hspace{5.6truecm}\displaystyle{\left.``\,D_{\alpha}\subset D_{\lambda}\,\,{\rm and}\,\,
\vert\lambda_x(v)\vert>\vert\alpha_{\hat x}(v)\vert\,",\,\,\,\,j\in{\Bbb Z}\right\}}\\
\displaystyle{
\bigcup\left\{\left.\frac{{\rm arctanh}(\lambda_x(v)/\alpha_{\hat x}(v))+\left(j+\frac 12\right)\pi
\sqrt{-1}}{\alpha_{\hat x}(v)}\,\right\vert\,(\alpha,\lambda)\in{\cal R}_M\times{\cal A}_M\,\,
{\rm s.t.}\,\,
\right.}\\
\hspace{5.6truecm}\displaystyle{\left.
``\,D_{\alpha}\subset D_{\lambda}\,\,{\rm and}\,\,
\vert\lambda_x(v)\vert<\vert\alpha_{\hat x}(v)\vert\,",\,\,\,\,j\in{\Bbb Z}\right\}.}
\end{array}$$
Hence it follows from (i) of Theorem A and (i) of Theorem B that $M$ has parallel focal structure 
(resp. complex equifocal) in the case where $G/K$ is of compact type (resp. non-compact type).  
Therefore the statements (ii) and (iii) of Theorem B follow from Remark 2.1 and Fact 2.3.  
\hspace{1.5truecm}q.e.d.

\newpage


\centerline{
\unitlength 0.1in
\begin{picture}( 48.6400, 23.1000)( -0.2700,-26.5000)
%
\special{pn 8}%
\special{ar 5098 2798 2330 1098  3.5314136 4.6003026}%
%
\special{pn 8}%
\special{ar 4866 2598 2604 1404  3.5945762 4.6225570}%
%
\special{pn 8}%
\special{ar 3172 1838 144 336  5.0352669 6.2831853}%
\special{ar 3172 1838 144 336  0.0000000 1.3298553}%
%
\special{pn 8}%
\special{ar 3230 1838 120 322  1.4934601 4.7123890}%
%
\special{pn 8}%
\special{ar 4148 1532 80 270  4.6196327 6.2831853}%
\special{ar 4148 1532 80 270  0.0000000 1.5152978}%
%
\special{pn 8}%
\special{ar 4148 1542 72 270  1.5707963 4.7123890}%
%
\special{pn 8}%
\special{pa 3462 1226}%
\special{pa 3462 2650}%
\special{fp}%
\special{pa 2004 340}%
\special{pa 2004 340}%
\special{fp}%
%
\special{pn 8}%
\special{pa 3462 2640}%
\special{pa 3100 2324}%
\special{fp}%
%
\special{pn 8}%
\special{pa 2890 2112}%
\special{pa 3020 2228}%
\special{fp}%
%
\special{pn 8}%
\special{pa 2890 846}%
\special{pa 3470 1226}%
\special{fp}%
%
\special{pn 8}%
\special{pa 2890 846}%
\special{pa 2890 1628}%
\special{fp}%
%
\special{pn 8}%
\special{pa 2890 1754}%
\special{pa 2890 2112}%
\special{fp}%
%
\special{pn 13}%
\special{pa 3260 1564}%
\special{pa 3334 1290}%
\special{fp}%
\special{sh 1}%
\special{pa 3334 1290}%
\special{pa 3298 1350}%
\special{pa 3320 1342}%
\special{pa 3336 1360}%
\special{pa 3334 1290}%
\special{fp}%
%
\special{pn 20}%
\special{sh 1}%
\special{ar 3260 1574 10 10 0  6.28318530717959E+0000}%
\special{sh 1}%
\special{ar 3260 1574 10 10 0  6.28318530717959E+0000}%
%
\special{pn 20}%
\special{sh 1}%
\special{ar 4154 1258 10 10 0  6.28318530717959E+0000}%
\special{sh 1}%
\special{ar 4154 1258 10 10 0  6.28318530717959E+0000}%
%
\special{pn 13}%
\special{pa 4154 1258}%
\special{pa 4148 994}%
\special{fp}%
\special{sh 1}%
\special{pa 4148 994}%
\special{pa 4130 1062}%
\special{pa 4148 1048}%
\special{pa 4170 1060}%
\special{pa 4148 994}%
\special{fp}%
%
\special{pn 8}%
\special{pa 4382 1046}%
\special{pa 4382 2250}%
\special{fp}%
%
\special{pn 8}%
\special{pa 4382 2250}%
\special{pa 4026 1860}%
\special{fp}%
%
\special{pn 8}%
\special{pa 3978 1796}%
\special{pa 3914 1722}%
\special{fp}%
%
\special{pn 8}%
\special{pa 3922 1742}%
\special{pa 3922 1354}%
\special{fp}%
%
\special{pn 8}%
\special{pa 3922 1248}%
\special{pa 3922 658}%
\special{fp}%
%
\special{pn 8}%
\special{pa 3922 658}%
\special{pa 4382 1058}%
\special{fp}%
%
\special{pn 8}%
\special{ar 1738 2704 1808 1086  3.5964914 3.6047901}%
\special{ar 1738 2704 1808 1086  3.6296864 3.6379851}%
\special{ar 1738 2704 1808 1086  3.6628814 3.6711802}%
\special{ar 1738 2704 1808 1086  3.6960764 3.7043752}%
\special{ar 1738 2704 1808 1086  3.7292715 3.7375702}%
\special{ar 1738 2704 1808 1086  3.7624665 3.7707652}%
\special{ar 1738 2704 1808 1086  3.7956615 3.8039603}%
\special{ar 1738 2704 1808 1086  3.8288565 3.8371553}%
\special{ar 1738 2704 1808 1086  3.8620515 3.8703503}%
\special{ar 1738 2704 1808 1086  3.8952466 3.9035453}%
\special{ar 1738 2704 1808 1086  3.9284416 3.9367403}%
\special{ar 1738 2704 1808 1086  3.9616366 3.9699354}%
\special{ar 1738 2704 1808 1086  3.9948316 4.0031304}%
\special{ar 1738 2704 1808 1086  4.0280266 4.0363254}%
\special{ar 1738 2704 1808 1086  4.0612217 4.0695204}%
\special{ar 1738 2704 1808 1086  4.0944167 4.1027154}%
\special{ar 1738 2704 1808 1086  4.1276117 4.1359105}%
\special{ar 1738 2704 1808 1086  4.1608067 4.1691055}%
\special{ar 1738 2704 1808 1086  4.1940017 4.2023005}%
\special{ar 1738 2704 1808 1086  4.2271968 4.2354955}%
\special{ar 1738 2704 1808 1086  4.2603918 4.2686905}%
\special{ar 1738 2704 1808 1086  4.2935868 4.3018856}%
\special{ar 1738 2704 1808 1086  4.3267818 4.3350806}%
\special{ar 1738 2704 1808 1086  4.3599768 4.3682756}%
\special{ar 1738 2704 1808 1086  4.3931719 4.4014706}%
\special{ar 1738 2704 1808 1086  4.4263669 4.4346656}%
\special{ar 1738 2704 1808 1086  4.4595619 4.4678607}%
\special{ar 1738 2704 1808 1086  4.4927569 4.5010557}%
\special{ar 1738 2704 1808 1086  4.5259520 4.5342507}%
\special{ar 1738 2704 1808 1086  4.5591470 4.5674457}%
\special{ar 1738 2704 1808 1086  4.5923420 4.6006408}%
%
\special{pn 8}%
\special{ar 1882 2608 1654 760  3.3671731 4.6103718}%
%
\special{pn 8}%
\special{ar 1916 2714 1500 632  3.2666352 3.2778975}%
\special{ar 1916 2714 1500 632  3.3116844 3.3229468}%
\special{ar 1916 2714 1500 632  3.3567337 3.3679960}%
\special{ar 1916 2714 1500 632  3.4017830 3.4130453}%
\special{ar 1916 2714 1500 632  3.4468323 3.4580946}%
\special{ar 1916 2714 1500 632  3.4918815 3.5031439}%
\special{ar 1916 2714 1500 632  3.5369308 3.5481931}%
\special{ar 1916 2714 1500 632  3.5819801 3.5932424}%
\special{ar 1916 2714 1500 632  3.6270294 3.6382917}%
\special{ar 1916 2714 1500 632  3.6720786 3.6833409}%
\special{ar 1916 2714 1500 632  3.7171279 3.7283902}%
\special{ar 1916 2714 1500 632  3.7621772 3.7734395}%
\special{ar 1916 2714 1500 632  3.8072264 3.8184888}%
\special{ar 1916 2714 1500 632  3.8522757 3.8635380}%
\special{ar 1916 2714 1500 632  3.8973250 3.9085873}%
\special{ar 1916 2714 1500 632  3.9423743 3.9536366}%
\special{ar 1916 2714 1500 632  3.9874235 3.9986858}%
\special{ar 1916 2714 1500 632  4.0324728 4.0437351}%
\special{ar 1916 2714 1500 632  4.0775221 4.0887844}%
\special{ar 1916 2714 1500 632  4.1225713 4.1338337}%
\special{ar 1916 2714 1500 632  4.1676206 4.1788829}%
\special{ar 1916 2714 1500 632  4.2126699 4.2239322}%
\special{ar 1916 2714 1500 632  4.2577192 4.2689815}%
\special{ar 1916 2714 1500 632  4.3027684 4.3140308}%
\special{ar 1916 2714 1500 632  4.3478177 4.3590800}%
\special{ar 1916 2714 1500 632  4.3928670 4.4041293}%
\special{ar 1916 2714 1500 632  4.4379163 4.4491786}%
\special{ar 1916 2714 1500 632  4.4829655 4.4942278}%
\special{ar 1916 2714 1500 632  4.5280148 4.5392771}%
\special{ar 1916 2714 1500 632  4.5730641 4.5843264}%
\special{ar 1916 2714 1500 632  4.6181133 4.6293757}%
\special{ar 1916 2714 1500 632  4.6631626 4.6711487}%
%
\special{pn 8}%
\special{ar 182 2736 1332 1604  5.5574234 5.8930628}%
%
\special{pn 8}%
\special{ar 118 2778 708 1150  5.3912904 5.8713180}%
%
\special{pn 20}%
\special{sh 1}%
\special{ar 682 2102 10 10 0  6.28318530717959E+0000}%
\special{sh 1}%
\special{ar 682 2102 10 10 0  6.28318530717959E+0000}%
%
\special{pn 20}%
\special{sh 1}%
\special{ar 1320 1912 10 10 0  6.28318530717959E+0000}%
\special{sh 1}%
\special{ar 1320 1912 10 10 0  6.28318530717959E+0000}%
%
\special{pn 13}%
\special{pa 682 2092}%
\special{pa 682 1722}%
\special{fp}%
\special{sh 1}%
\special{pa 682 1722}%
\special{pa 662 1790}%
\special{pa 682 1776}%
\special{pa 702 1790}%
\special{pa 682 1722}%
\special{fp}%
%
\special{pn 13}%
\special{pa 1320 1902}%
\special{pa 1310 1522}%
\special{fp}%
\special{sh 1}%
\special{pa 1310 1522}%
\special{pa 1292 1590}%
\special{pa 1312 1576}%
\special{pa 1332 1588}%
\special{pa 1310 1522}%
\special{fp}%
%
\special{pn 13}%
\special{pa 682 2102}%
\special{pa 790 2280}%
\special{fp}%
\special{sh 1}%
\special{pa 790 2280}%
\special{pa 774 2214}%
\special{pa 762 2234}%
\special{pa 738 2234}%
\special{pa 790 2280}%
\special{fp}%
%
\special{pn 13}%
\special{pa 1328 1922}%
\special{pa 1430 2100}%
\special{fp}%
\special{sh 1}%
\special{pa 1430 2100}%
\special{pa 1414 2032}%
\special{pa 1404 2054}%
\special{pa 1380 2052}%
\special{pa 1430 2100}%
\special{fp}%
%
\special{pn 8}%
\special{ar 2826 2124 1456 866  3.5271211 3.5374659}%
\special{ar 2826 2124 1456 866  3.5685004 3.5788452}%
\special{ar 2826 2124 1456 866  3.6098797 3.6202246}%
\special{ar 2826 2124 1456 866  3.6512590 3.6616039}%
\special{ar 2826 2124 1456 866  3.6926384 3.7029832}%
\special{ar 2826 2124 1456 866  3.7340177 3.7443625}%
\special{ar 2826 2124 1456 866  3.7753970 3.7857418}%
\special{ar 2826 2124 1456 866  3.8167763 3.8271211}%
\special{ar 2826 2124 1456 866  3.8581556 3.8685004}%
\special{ar 2826 2124 1456 866  3.8995349 3.9098797}%
\special{ar 2826 2124 1456 866  3.9409142 3.9512590}%
\special{ar 2826 2124 1456 866  3.9822935 3.9926384}%
\special{ar 2826 2124 1456 866  4.0236728 4.0340177}%
\special{ar 2826 2124 1456 866  4.0650521 4.0753970}%
\special{ar 2826 2124 1456 866  4.1064315 4.1167763}%
\special{ar 2826 2124 1456 866  4.1478108 4.1581556}%
\special{ar 2826 2124 1456 866  4.1891901 4.1995349}%
\special{ar 2826 2124 1456 866  4.2305694 4.2409142}%
\special{ar 2826 2124 1456 866  4.2719487 4.2822935}%
\special{ar 2826 2124 1456 866  4.3133280 4.3236728}%
\special{ar 2826 2124 1456 866  4.3547073 4.3650521}%
%
\special{pn 8}%
\special{pa 2334 1322}%
\special{pa 2406 1300}%
\special{dt 0.045}%
\special{sh 1}%
\special{pa 2406 1300}%
\special{pa 2336 1300}%
\special{pa 2356 1316}%
\special{pa 2348 1338}%
\special{pa 2406 1300}%
\special{fp}%
%
\special{pn 8}%
\special{pa 538 1510}%
\special{pa 610 1944}%
\special{dt 0.045}%
\special{sh 1}%
\special{pa 610 1944}%
\special{pa 618 1876}%
\special{pa 600 1892}%
\special{pa 578 1882}%
\special{pa 610 1944}%
\special{fp}%
%
\special{pn 8}%
\special{pa 1198 1278}%
\special{pa 1238 1764}%
\special{dt 0.045}%
\special{sh 1}%
\special{pa 1238 1764}%
\special{pa 1252 1696}%
\special{pa 1234 1712}%
\special{pa 1212 1700}%
\special{pa 1238 1764}%
\special{fp}%
\put(4.6300,-14.5800){\makebox(0,0)[lb]{$\Sigma^{\alpha}_x$}}%
\put(11.3300,-12.4800){\makebox(0,0)[lb]{$\Sigma^{\alpha}_y$}}%
\put(6.4900,-16.6900){\makebox(0,0)[lb]{$\widetilde w_x$}}%
\put(12.5400,-14.5800){\makebox(0,0)[lb]{$\widetilde w_y$}}%
\put(8.1900,-23.4400){\makebox(0,0)[lt]{$(\xi_{\alpha})_{\hat x}$}}%
\put(14.1500,-21.7600){\makebox(0,0)[lt]{$(\xi_{\alpha})_{\hat y}$}}%
\put(38.3300,-9.9400){\makebox(0,0)[rb]{$(\xi_{\alpha})_{\hat y}$}}%
\put(26.8000,-10.8900){\makebox(0,0)[rb]{$(\xi_{\alpha})_{\hat x}$}}%
\put(27.4500,-8.5700){\makebox(0,0)[rb]{$\Sigma^{\alpha}_x$}}%
\put(38.1700,-7.5100){\makebox(0,0)[rb]{$\Sigma^{\alpha}_y$}}%
%
\special{pn 8}%
\special{pa 2714 1046}%
\special{pa 3294 1438}%
\special{dt 0.045}%
\special{sh 1}%
\special{pa 3294 1438}%
\special{pa 3250 1384}%
\special{pa 3250 1408}%
\special{pa 3228 1416}%
\special{pa 3294 1438}%
\special{fp}%
%
\special{pn 8}%
\special{pa 3850 974}%
\special{pa 4148 1164}%
\special{dt 0.045}%
\special{sh 1}%
\special{pa 4148 1164}%
\special{pa 4102 1110}%
\special{pa 4102 1134}%
\special{pa 4080 1144}%
\special{pa 4148 1164}%
\special{fp}%
%
\special{pn 8}%
\special{pa 3834 710}%
\special{pa 4004 884}%
\special{dt 0.045}%
\special{sh 1}%
\special{pa 4004 884}%
\special{pa 3972 822}%
\special{pa 3966 846}%
\special{pa 3942 850}%
\special{pa 4004 884}%
\special{fp}%
%
\special{pn 8}%
\special{pa 2770 826}%
\special{pa 3010 1100}%
\special{dt 0.045}%
\special{sh 1}%
\special{pa 3010 1100}%
\special{pa 2982 1036}%
\special{pa 2976 1060}%
\special{pa 2952 1062}%
\special{pa 3010 1100}%
\special{fp}%
\put(33.0800,-15.4200){\makebox(0,0)[lt]{$x$}}%
\put(41.8700,-12.2600){\makebox(0,0)[lb]{$y$}}%
%
\special{pn 8}%
\special{pa 3696 2314}%
\special{pa 3316 1932}%
\special{dt 0.045}%
\special{sh 1}%
\special{pa 3316 1932}%
\special{pa 3350 1994}%
\special{pa 3354 1970}%
\special{pa 3378 1966}%
\special{pa 3316 1932}%
\special{fp}%
%
\special{pn 8}%
\special{pa 4600 1986}%
\special{pa 4228 1606}%
\special{dt 0.045}%
\special{sh 1}%
\special{pa 4228 1606}%
\special{pa 4260 1668}%
\special{pa 4264 1644}%
\special{pa 4288 1640}%
\special{pa 4228 1606}%
\special{fp}%
\put(46.3900,-19.8600){\makebox(0,0)[lt]{$L^{\alpha,2\alpha}_y$}}%
\put(37.1100,-23.4400){\makebox(0,0)[lt]{$L^{\alpha,2\alpha}_x$}}%
\put(21.4000,-13.1000){\makebox(0,0)[rb]{in fact}}%
%
\special{pn 8}%
\special{pa 4642 784}%
\special{pa 4542 1362}%
\special{dt 0.045}%
\special{sh 1}%
\special{pa 4542 1362}%
\special{pa 4574 1300}%
\special{pa 4552 1308}%
\special{pa 4534 1292}%
\special{pa 4542 1362}%
\special{fp}%
%
\special{pn 8}%
\special{pa 250 1670}%
\special{pa 466 2218}%
\special{dt 0.045}%
\special{sh 1}%
\special{pa 466 2218}%
\special{pa 460 2148}%
\special{pa 446 2168}%
\special{pa 424 2162}%
\special{pa 466 2218}%
\special{fp}%
\put(2.4300,-16.1900){\makebox(0,0)[rb]{$M$}}%
\put(46.4100,-7.4200){\makebox(0,0)[lb]{$M$}}%
%
\special{pn 4}%
\special{pa 130 2248}%
\special{pa 276 2070}%
\special{fp}%
\special{pa 154 2280}%
\special{pa 394 1990}%
\special{fp}%
\special{pa 180 2312}%
\special{pa 496 1928}%
\special{fp}%
\special{pa 204 2344}%
\special{pa 590 1878}%
\special{fp}%
\special{pa 228 2376}%
\special{pa 676 1836}%
\special{fp}%
\special{pa 254 2408}%
\special{pa 750 1808}%
\special{fp}%
\special{pa 278 2440}%
\special{pa 824 1780}%
\special{fp}%
\special{pa 304 2472}%
\special{pa 896 1756}%
\special{fp}%
\special{pa 328 2504}%
\special{pa 968 1730}%
\special{fp}%
\special{pa 354 2534}%
\special{pa 1030 1716}%
\special{fp}%
\special{pa 378 2566}%
\special{pa 1090 1706}%
\special{fp}%
\special{pa 404 2598}%
\special{pa 470 2516}%
\special{fp}%
\special{pa 534 2440}%
\special{pa 1150 1696}%
\special{fp}%
\special{pa 650 2362}%
\special{pa 1210 1686}%
\special{fp}%
\special{pa 746 2308}%
\special{pa 1268 1676}%
\special{fp}%
\special{pa 836 2260}%
\special{pa 1328 1668}%
\special{fp}%
\special{pa 918 2224}%
\special{pa 1386 1660}%
\special{fp}%
\special{pa 982 2208}%
\special{pa 1440 1654}%
\special{fp}%
\special{pa 1048 2190}%
\special{pa 1494 1650}%
\special{fp}%
\special{pa 1114 2170}%
\special{pa 1550 1646}%
\special{fp}%
\special{pa 1180 2154}%
\special{pa 1606 1640}%
\special{fp}%
\special{pa 1240 2144}%
\special{pa 1626 1678}%
\special{fp}%
\special{pa 1300 2134}%
\special{pa 1646 1716}%
\special{fp}%
\special{pa 1358 2124}%
\special{pa 1666 1754}%
\special{fp}%
\special{pa 1418 2114}%
\special{pa 1686 1790}%
\special{fp}%
\special{pa 1474 2108}%
\special{pa 1706 1828}%
\special{fp}%
\special{pa 1530 2102}%
\special{pa 1726 1866}%
\special{fp}%
\special{pa 1586 2096}%
\special{pa 1746 1904}%
\special{fp}%
\special{pa 1642 2092}%
\special{pa 1766 1940}%
\special{fp}%
\special{pa 1696 2088}%
\special{pa 1786 1980}%
\special{fp}%
%
\special{pn 4}%
\special{pa 1750 2084}%
\special{pa 1806 2016}%
\special{fp}%
\special{pa 1808 2078}%
\special{pa 1826 2054}%
\special{fp}%
\put(10.5000,-25.7000){\makebox(0,0)[lt]{($\hat x\in\pi^{-1}(x)$)}}%
\put(16.2000,-23.8000){\makebox(0,0)[lt]{($\hat y\in\pi^{-1}(y)$)}}%
\end{picture}%
\hspace{0truecm}}

\vspace{0.8truecm}

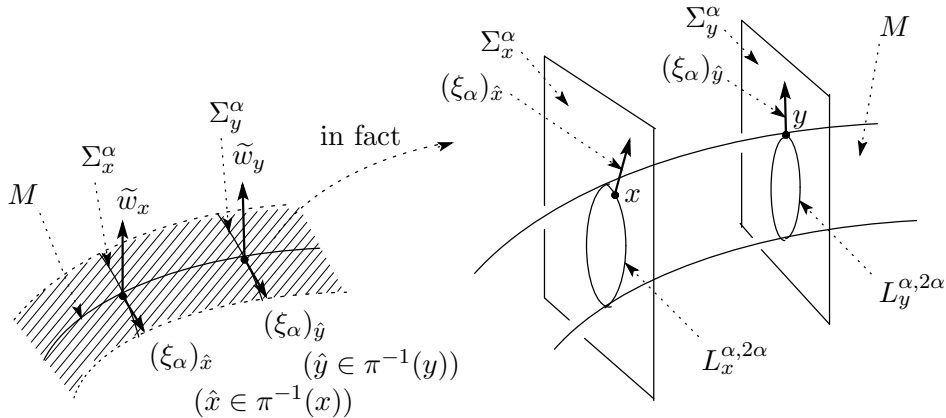
\centerline{\bf Figure 1.}

\vspace{0.5truecm}

Next we shall prove Theorem C.  

\vspace{0.5truecm}

\noindent
{\it Proof of Theorem C.} According to the proof of (ii) of Theorem B, $M$ is equifocal.  
Hence, according to Fact 2.1, it follows from the additional conditions in Theorem C 
that $M$ is a principal orbit of a Hermann action on $G/K$.  
Denote by $H$ this Hermann action.  
Furthermore, it follows from the additional assumption in Theorem B that 
the $H$-action satisfies the condition ${\triangle'}^V_+\cap{\triangle'}^H_+=\emptyset$, 
where ${\triangle'}^V_+$ and ${\triangle'}^H_+$ are the systems determined by the triple 
$(H,G,K)$ as in [21].  Hence, according to Proposition 4.39 in [13], the $H$-action is 
conjugate to the isotropy action of $G/K$.  That is, $M$ is congruent to a principal orbit 
of the isotropy action.  
\hspace{1.5truecm}q.e.d.

\vspace{0.5truecm}

Finally we shall prove Theorem D.  

\vspace{0.5truecm}

\noindent
{\it Proof of Theorem D.} According to (iii) of Theorem B, $M$ is isoparametric.  
According to the result in [23], it follows from this fact and the additional assumptions 
in Theorem D that $M$ is a principal orbit of a Hermann action on $G/K$.  \hspace{1.5truecm}q.e.d.

\section{Examples} 
Let $G/K$ be a symmetric space of non-compact type such that 
the multiplicity of each root of the (restricted) root system of the symmetric pair $(G,K)$ 
is greater than one.  Then the principal orbits of the isotropy actions of $G/K$ are 
complete curvature-adapted submanifolds as in Theorem B.  
Also, the principal orbits of Hermann actions $H\curvearrowright G/K$'s (of cohomogeneity 
two) in Table 1 and their dual actions $H^d\curvearrowright G^d/K$'s are complete curvature-adapted 
submanifolds as in Theorem B (see [8,18,25]), where $G^d/K$ is the compact dual of $G/K$.  
Note that the dual actions are conjugate to the isotropy actions of $G^d/K$ by Proposition 4.39 of [13].  
The Hermann actions in Table 1 are all of Hermann actions 
on irreducible rank two symmetric spaces of non-compact type such that their principal 
orbits satisfy all the hypothesis in Theorem B and that they are not conjugate to the 
isotropy actions (see [21,24]).  

\vspace{0.5truecm}

$$\begin{tabular}{|c|c|c|}
\hline
{\scriptsize $H\curvearrowright G/K$} & {\scriptsize ${\rm dim}\,M$} & 
{\scriptsize $\sharp\,{\cal A}_M$}\\
\hline
{\scriptsize$Sp(1,2)\curvearrowright SU^{\ast}(6)/Sp(3)$} & {\scriptsize $12$} & {\scriptsize $3$}\\
\hline
{\scriptsize$SO_0(2,3)\curvearrowright SO(5,{\Bbb C})/SO(5)$} & {\scriptsize $8$} & 
{\scriptsize $4$}\\
\hline
{\scriptsize$Sp(2,{\Bbb R})\curvearrowright Sp(2,{\Bbb C})/Sp(2)$} & {\scriptsize $8$} & 
{\scriptsize $4$}\\
\hline
{\scriptsize$Sp(1,1)\curvearrowright Sp(2,{\Bbb C})/Sp(2)$} & {\scriptsize $8$} & 
{\scriptsize $4$}\\
\hline
{\scriptsize$F_4^{-20}\curvearrowright E_6^{-26}/F_4$} & {\scriptsize $24$} & {\scriptsize $3$}\\
\hline
{\scriptsize$G_2^2\curvearrowright G_2^{\bf C}/G_2$} & {\scriptsize $12$} & {\scriptsize $6$}\\
\hline
\end{tabular}$$

\vspace{0.1truecm}

\centerline{($M\,\,:\,\,$a principal orbit of $H\curvearrowright G/K$)}

\vspace{0.2truecm}

\centerline{{\bf Table 1.}}

\vspace{0.8truecm}

\noindent
{\bf References}

\vspace{0.4truecm}


%
%

\noindent
[1] J. Berndt, Real hypersurfaces with constant principal curvatures 
in complex 

hyperbolic space, J. Reine Angew. Math. {\bf 395} 132-141 (1989). 

\noindent
[2] J. Berndt, Real hypersurfaces in quaternionic space forms, 
J. Reine Angew. 

Math. {\bf 419} 9-26 (1991). 

\noindent
[3] J. Berndt and L. Vanhecke, 
Curvature adapted submanifolds, 
Nihonkai Math. 

J. {\bf 3} 177-185 (1992).

\noindent
[4] R.A. Blumenthal and J.J. Hebda, Complementary distributions which preserve 

the leaf geometry and applications to totally geodesic foliations, Quart. J. 

Math. Oxford. Ser. (2) {\bf 35} 383-392 (1984).  

\noindent
[5] T.E. Cecil and P.J. Ryan, Focal sets and real hypersurfaces in complex proje-

ctive space, Trans. Amer. Math. Soc. {\bf 269} 481-499 (1982).

\noindent
[6] U. Christ, Homogeneity of equifocal submanifolds, J. Differential Geometry 

{\bf 62} 1-15 (2002).

\noindent
[7] P.B. Eberlein, Geometry of nonpositively curved manifolds, Chicago Lectu-

res in Mathematics, Univerity of Chicago Press, Chicago, IL, 1996.

\noindent
[8] O. Goertsches and G. Thorbergsson, On the Geometry of the orbits of 

Hermann actions, Geom. Dedicata {\bf 129} 101-118 (2007).

\noindent
[9] C. Gorodski and E. Heintze, 
Homogeneous structures and rigidity of isopara-

metric submanifolds in Hilbert space, J. Fixed Point Theory Appl. {\bf 11} (2012) 

93-136.

\noindent
[10] E. Heintze, X. Liu and C. Olmos, Isoparametric submanifolds and a Cheva-

lley type restriction theorem, Integrable systems, geometry, and topology, 

151-190, AMS/IP Stud. Adv. Math. 36, Amer. Math. Soc., Providence, RI, 

2006.

\noindent
[11] E. Heintze, R.S. Palais, C.L. Terng and G. Thorbergsson, Hyperpolar acti-

ons on symmetric spaces, Geometry, topology and physics for Raoul Bott 

(ed. S. T. Yau), Conf. Proc. Lecture Notes Geom. Topology {\bf 4}, Internat. 

Press, Cambridge, MA, 1995 pp214-245.

\noindent
[12] S. Helgason, 
Differential geometry, Lie groups and symmetric spaces, Acad-

emic Press, New York, 1978.

\noindent
[13] O. Ikawa, The geometry of symmetric triad and orbit spaces of Hermann 

actions, J. Math. Soc. Japan {\bf 63} 79-139 (2011).  

\noindent
[14] N. Innami, Y. Mashiko and K. Shiohama, Metric spheres in the projective 

spaces with constant holomorphic sectional curvautre, Tsukuba J. Math. {\bf 35} 

79-90 (2011).  

\noindent
[15] N. Innami, Y. Itokawa and K. Shiohama, Complete real hypersurfaces and 

special ${\Bbb K}$-line bundles in ${\Bbb K}$-hyperbolic spaces, Intern. J. 
Math. {\bf 24} 1350082, 

14 pp. (2013).  

\noindent
[16] N. Koike, 
Submanifold geometries in a symmetric space of non-compact 

type and a pseudo-Hilbert space, Kyushu J. Math. {\bf 58} 167-202 (2004).

\noindent
[17] N. Koike, 
Complex equifocal submanifolds and infinite dimensional anti-

Kaehlerian isoparametric submanifolds, Tokyo J. Math. {\bf 28} 201-247 
(2005).

\noindent
[18] N. Koike, Actions of Hermann type and proper complex equifocal subma-

nifolds, Osaka J. Math. {\bf 42} 599-611 (2005).  

\noindent
[19] N. Koike, 
On curvature-adapted and proper complex equifocal submanifo-

lds, Kyungpook Math. J. {\bf 50} 509-536 (2010).

\noindent
[20] N. Koike, 
Examples of a complex hyperpolar action without singular orbit, 

Cubo A Math. J. {\bf 12} 131-147 (2010).

\noindent
[21] N. Koike, 
Collapse of the mean curvature flow for equifocal submanifolds, 

Asian J. Math. {\bf 15} 101-128. (2011)

\noindent
[22] N. Koike, 
An Cartan type identity for isoparametric hypersurfaces in sym-

metric spaces, Tohoku Math. J. (to appear) (arXiv:math. DG/1010.1652

v3).

\noindent
[23] N. Koike, 
The classifications of certain kind of isoparametric submanifolds 

in non-compact symmetric spaces, arXiv:math.DG/1209.1933v1.

\noindent
[24] A. Kollross, A Classification of hyperpolar and cohomogeneity one actions, 

Trans. Amer. Math. Soc. {\bf 354} 571-612 (2001).

\noindent
[25] A. Kollross, Duality of symmetric spaces and polar actions, J. Lie Theory 

{\bf 21} 961-986 (2011).

\noindent
[26] S. Montiel, Real hypersurfaces of a complex hyperbolic space, J. Math. Soc. 

Japan {\bf 37} 515-535 (1985).

\noindent
[27] M. Ortega and J.D. P$\acute e$rez, On the Ricci tensor of a real hypersurface of 

quaternionic hyperbolic space, manuscripta math. {\bf 93} 49-57 (1997).

\noindent
[28] T. Murphy, Curvature-adapted submanifolds of symmetric spaces, Indiana 

Univ. Math. J. {\bf 61} 831-847 (2012). 

\noindent
[29] R.S. Palais and C.L. Terng, Critical point theory and submanifold geometry, 

Lecture Notes in Math. {\bf 1353}, Springer, Berlin, 1988.

\noindent
[30] T. Sakai, Riemannian Geometry, Translations of Mathematical Monographs, 

149. American Math. Soc., Providence, RI,1996.

\noindent
[31] C.L. Terng, Isoparametric submanifolds and their Coxeter groups, J. Differe-

ntial Geometry {\bf 21} 79-107 (1985).

\noindent
[32] C.L. Terng and G. Thorbergsson, Submanifold geometry in symmetric spaces, 

J. Differential Geometry {\bf 42} 665-718 (1995).


\vspace{0.5truecm}

\rightline{Department of Mathematics, Faculty of Science, }
\rightline{Tokyo University of Science}
\rightline{1-3 Kagurazaka Shinjuku-ku,}
\rightline{Tokyo 162-8601, Japan}
\rightline{(e-mail: koike@ma.kagu.tus.ac.jp)}

\end{document}